\documentclass{amsart}
\usepackage{amsmath,amsthm}
\usepackage{amsfonts,amssymb}

\newtheorem{theorem}{Theorem}[section]
\newtheorem{corollary}[theorem]{Corollary}

\newtheorem{proposition}[theorem]{Proposition}
\newtheorem{definition}[theorem]{Definition}

\theoremstyle{remark}

 \def\ab{{\mathbf a}}
 
 \def\cb{{\mathbf c}}
 \def\nb{{\mathbf n}}

 \def\xb{{\mathbf x}}

 \def\CL{{\mathcal L}}
 
 \def\CP{{\mathcal P}}     
 \def\CS{{\mathcal S}}     
 \def\CV{{\mathcal V}}

 \def\NN{{\mathbb N}}
 \def\PP{{\mathbb P}}
 
 \def\RR{{\mathbb R}}

        \def\sspan{\operatorname{span}}
        
        \def\dim{\operatorname{dim}}
        \def\rank{\operatorname{rank}}
        \def\diag{\operatorname{diag}}
        
        \def\LT{\operatorname{\sf LT}}

\newcommand{\wt}{\widetilde}

\input epsf.tex
 
\begin{document}

\title{On discrete orthogonal polynomials of several variables}  
\author{Yuan Xu}
\address{Department of Mathematics\\ University of Oregon\\
    Eugene, Oregon 97403-1222.}\email{yuan@math.uoregon.edu}

\date{\today}
\thanks{Work partially supported by the National Science Foundation under 
Grant DMS-0201669}
                         
\begin{abstract}
Let $V$ be a set of isolated points in $\RR^d$. Define a linear functional 
$\CL$ on the space of real polynomials restricted on $V$, $\CL f = 
\sum_{x \in V} f(x)\rho(x)$, where $\rho$ is a nonzero function on $V$. 
Polynomial subspaces that contain discrete orthogonal polynomials with 
respect to the bilinear form $\langle f,g\rangle = \CL(f g)$ are 
identified. One result shows that the discrete orthogonal polynomials still 
satisfy a three-term relation and Favard's theorem holds in this general
setting. 
\end{abstract}

\maketitle                      
 
\section{Introduction}
\setcounter{equation}{0}

Discrete orthogonal polynomials appear naturally in combinatorics, genetics,
statistics and various areas in applied mathematics (see, for example, 
\cite{KM,NSU}). In one variable they have been studied extensively. 
Let $V$ be a set of isolated points on the real line, its cardinality $|V|$ 
is either finite or countable. Let $w$ be a real positive function on $V$. 
With respect to the bilinear form $\langle f,g \rangle = \sum_{x\in V} f(x)g(x)
 w(x)$, there is a sequence of orthogonal polynomials $\{p_n:0\le n \le |V|\}$
on $V$ with $\langle p_n, p_m \rangle =0$ for $n\ne m$ (for example, using 
Gram-Schmidt process). These are the discrete orthogonal polynomials. Their
structure is similar to that of the usual continuous orthogonal polynomials. 
For example, every sequence of discrete orthogonal polynomials satisfies a 
three-term relation, 
\begin{equation}
  x p_n = a_n p_{n+1} + b_n p_n + c_n p_{n-1}, \qquad  0 \le n \le |V|-1,
\end{equation}
where $a_n$, $b_n$ and $c_n$ are real numbers. Furthermore, according to 
Favard's theorem, the three-term relation essentially characterizes the 
orthogonality of polynomials. 

Discrete orthogonal polynomials of several variables are far less studied. 
Their orthogonal structure is much more complicated than that of one variable. 
Even some basic problems have not been addressed. Let us first fix some 
notation. We use the standard multiindex notation: for $x \in \RR^d$ and 
$\alpha\in \NN_0^d$, write $x^\alpha=x_1^{\alpha_1}\cdots 
x_d^{\alpha_d}$. This is a monomial of (total) degree $|\alpha|:= \alpha_1+ 
\ldots +\alpha_d$. Let $\Pi^d = \RR[x_1,\ldots,x_d]$ and let $\Pi_n^d$ be the 
subspace of polynomials of degree at most $n$. Denote by $\CP_n^d$ the space 
of homogeneous polynomials of degree $n$. It is well known that 
$$
  \dim \CP_n^d = \binom{n+d-1}{d-1} \quad \hbox{and} \quad
  \dim \Pi_n^d = \binom{n+d}{d}.
$$
Let $\CL$ be a linear functional defined on $\Pi^d$, such that a basis of 
orthogonal polynomials $\{P_\alpha:|\alpha| = n, \alpha \in \NN_0^d, n 
\ge 0\}$, where $P_\alpha \in \Pi_n^d$, exists with respect to the bilinear
form $\langle f, g\rangle = \CL (f g)$,
in the sense that $\langle P_\alpha, P_\beta \rangle = 0$ if $|\alpha|\ne 
|\beta|$. Let $\CV_n^d = \sspan\{P_\alpha:|\alpha| = n\}$ be the space of 
orthogonal polynomials of total degree $n$. Then $\dim \CV_n^d =\dim \CP_n^d$.
In this case, if we adopt the point of view that the orthogonality holds in 
terms of the subspaces $\CV_n^d$, not in terms of particular bases of 
$\CV_n^d$, then we can have an analog of a three-term relation. Let 
$\PP_n = \{P_\alpha: |\alpha|=n\}$; we also use $\PP_n$ to denote a column 
vector, in which the elements are ordered according to a fixed monomial order.
Then the following three-term relation holds,
\begin{equation}\label{eq:1.2}
  x_i \PP_n = A_{n,i} \PP_{n+1} + B_{n,i} \PP_n + C_{n,i} \PP_{n-1}, \qquad 
     1 \le i \le d, \quad n \ge 0,
\end{equation}
where $A_{n,i}$, $B_{n,i}$ and $C_{n,i}$ are matrices of appropriate
dimensions, 
and $\PP_{-1}:=0$. Furthermore, there is an analogue of Favard's theorem 
(\cite{X93}). For the general theory of orthogonal polynomials of several 
variables, we refer to Chapter 3 of \cite{DX}. 

Let $V$ be a set of isolated points in $\RR^d$. Again we denote by $|V|$ the
cardinality of $V$, which can be finite or countable. The orthogonal 
polynomials on $V$ depends on the structure of the polynomial ideal $I(V)$ 
that has $V$ as its variety, 
$$
   I(V) = \{p \in \RR[x_1,\ldots,x_d]: p(x) =0, \forall x \in V\}.
$$ 
The discrete orthogonal polynomials on $V$ can only consist of polynomials
that do not belong to $I(V)$. Let $\CL$ be defined by $\CL f = 
\sum_{x\in V} f(x) W(x)$, where $W$ is a real function on $V$, $W(x) \ne
0$ for all $x \in V$, and $ \sum_{x \in V} |x^\alpha| |W(x)| < \infty$
for all  $\alpha \in \NN_0^d$
in the case where $V$ is a countable set. Only if $I(V) = \{0\}$, can a 
complete basis $\{P_\alpha: \alpha \in \NN_0^d\}$ of $\Pi^d$ 
with respect to $\CL$ exist and the discussion in the previous paragraph 
apply. If the ideal $I(V)$ is non-trivial, for example, when $|V|$ is finite, 
then we need to understand the subspace $\RR[V] \cong \Pi^d/I(V)$ in order to 
define orthogonal polynomials. In the case $\RR[V] = \Pi_N^d$, for example,
little extra work is needed; the three-term relation in the form 
\eqref{eq:1.2} holds for $0 \le n \le N$. This is the case of straightforward 
extension of discrete orthogonal polynomials in one variable. However, even in 
the case where $V$ is a product of two point sets $X$ and $Y$ in one variable, 
$V = X\times Y$ with $|V| = N$, $\RR[V]$ consists of only a subspace of 
$\Pi_N^d$. In general, the space $\RR[V]$ can be rather complicated and care
is needed for the definition of orthogonal polynomials on a discrete set $V$. 

The purpose of the present study is to define discrete orthogonal polynomials 
in this general setting. The polynomial subspaces for which the discrete 
orthogonal polynomials exist are identified and the three-term relation and 
Favard's theorem are established. In the following section we 
discuss the structure of polynomial subspace on $V$. Discrete orthogonal 
polynomials on $V$  are studied in Section 3. Various examples will be 
given in the paper, some are given in terms of the classical discrete 
orthogonal polynomials, such as Hahn polynomials, in Section 4. It is our hope 
that this study can help to clarify some of the basic questions in the theory 
of discrete orthogonal polynomials. 

The author would like to thank Paul Terwilliger for the discussion on the 
Leonard pairs and tridiagonal pairs (see \cite{T2,T1} and the reference 
there), during the 7th International Symposium on Orthogonal Polynomials 
and Special Functions, held in Copenhagen in August, 2003, which suggests
this study. 

\section{Polynomial spaces on $V$}
\setcounter{equation}{0}

First we need to understand the structure of the quotient ideal $\Pi^d /I$,
where $I: = I(V)$ and $V$ is a set of isolated points in $\RR^d$, finite or 
countable. Most of the results below also hold if $V$ has finitely many 
accumulation points. We review some results about ideals and varieties, our 
basic reference is \cite{CLO}. 

For $f, g \in \Pi^d$, we say that $f$ is congruent to $g$ modulo $I$, 
written as $f = g \mod I$, if and only if $f-g \in I$. If $|V|$ is finite, 
then it is known that the codimension of $I$ is equal to $|V|$; that is, 
$\dim \Pi^d/I(V) =|V|$. Let $\RR[V]$ denote the collection of polynomial 
functions $\phi: V \mapsto \RR$. This is a commutative ring and it is 
isomorphic to the quotient ring $\Pi^d/I(V)$. We can identify $\RR[V]$ 
with $\Pi^d /I$ as there is an one-to-one correspondence between $\phi \in 
\RR[V]$ and $[\phi] = \{g \in \Pi^d: g=\phi  \mod I\}$. It is possible to 
say more about this space. For a fixed monomial order, we denote by $\LT(f)$ 
the leading monomial term for any polynomial $f \in \Pi^d$; that is, if $f = 
\sum c_\alpha x^\alpha$, then $\LT(f) = c_\beta x^\beta$, where $x^\beta$ 
is the leading monomial among all monomials appearing in $f$ and 
$c_\beta\ne 0$. 
For a polynomial ideal $I$ other than $\{1\}$, we denote by $\LT(I)$ the 
leading terms of $I$, that is,  
$$
   \LT(I) = \{ c \xb^\alpha \vert \hbox{there exists $f\in I$ with
        $\LT(f) = c \xb^\alpha$}\}.  
$$
We further denote by $\langle \LT(I) \rangle$ the ideal generated by the 
leading terms of $\LT(f)$ for all $f \in I \setminus \{0\}$. According to 
the Hilbert basis theorem, every polynomial ideal has a finite basis. A set 
$\{g_1, \ldots, g_t\}$ is called a Gr\"obner basis of $I$ if
$$
  \langle \LT(g_1), \ldots, \LT(g_t) \rangle = \langle \LT(I) \rangle. 
$$
It is known that every polynomial ideal has a Gr\"obner basis. Such a basis
enjoys many interesting properties that have important applications. For 
example, it is used to prove the following result (\cite[Chapt. 5]{CLO}).

\begin{proposition} \label{prop:2.1}
Fix a monomial ordering on $\Pi^d$ and let $I \subset \Pi^d$ be an ideal. Then 
there is an isomorphism between $\Pi^d / I$ and the space
$$
\CS_I := \sspan \{\xb^\alpha \vert \xb^\alpha \notin \langle \LT(I)\rangle\}. 
$$
More precisely, every $f \in \Pi^d$ is congruent modulo $I$ to a unique 
polynomial $r \in \CS_I$.  
\end{proposition}

In fact, the polynomial $r$ is the remainder of $f$ on division by $I$. For
$I = \langle f_1, \ldots, f_M \rangle$ and a fixed monomial order, the 
division algorithm states that for every $f \in \Pi^d$, there exist $p_i$ 
and $r$ in $\Pi^d$ such that $f = \sum p_i f_i + r$, where $r \in \CS_I$
and no term of $r$ is divisible by any of $\LT(f_1),\ldots,\LT(f_M)$. 
The remainder polynomial $r$ is unique if the basis $f_1, \ldots, f_M$ is a 
Gr\"obner basis. 

\begin{proposition} \label{prop:2.2}
Let $V$ be a set of isolated points in $\RR^d$ and $I = I(V)$. Let $\Lambda 
:=\Lambda(V)$ be the index set $\Lambda =\{\alpha:x^\alpha\notin \LT(I)\}$.
Then every polynomial $P \in \RR[V]$ can be written uniquely as 
$$
P(x) = \sum_{\alpha \in \Lambda} c_\alpha x^\alpha \mod I(V), 
  \qquad c_\alpha \in \RR,
$$ 
and the set $\Lambda$ satisfies the following property 
\begin{equation} \label{eq:2.1}
\alpha \in \Lambda \quad \hbox{implies} \quad \alpha-\beta\in \Lambda, \quad 
\hbox{whenever  $\alpha - \beta \in \NN_0^d$ and $\beta \in \NN_0^d$}.
\end{equation}
\end{proposition} 

\begin{proof}
For $I = I(V)$, we can take $\CS_I$ in Proposition \ref{prop:2.1} as $\RR[V]$,
modulus $I$ if needed. The definition shows that we can write $\CS_I = 
\sspan\{x^\alpha:\alpha \in \Lambda\}$. Hence, every polynomial $P$ in $\RR[V]$
has the stated representation. Since the ideal $\langle \LT(I)\rangle$ is a
monomial ideal, $x^\alpha \in \langle \LT(I)\rangle$ implies $x^{\alpha+\beta} 
\in \langle \LT(I)\rangle$ for any $\beta \in \NN_0^d$. Consequently, it 
follows that the set $\Lambda$ satisfies the property \eqref{eq:2.1}. 
\end{proof}

In the following, we shall drop modulus $I$ and use $\RR[V]$ to denote the
space 
\begin{equation} \label{eq:2.2}
\RR[V] = \sspan\{x^\alpha:\alpha \in \Lambda(V)\}.
\end{equation}
This abuse of notation should not cause problems.

We should point out that the set $\Lambda$ is not unique, since all equations
actually hold under congruence modulo $I$. In fact, Gr\"obner bases are not
unique, since the choice of monomial orders matters. There are in fact many 
different representations of elements in $\RR[V]$. What is of interest is the 
property \eqref{eq:2.1} satisfied by $\Lambda$. 

\medskip\noindent
{\bf Example 2.1} Consider the set $V = \{(0,0),(0,1),(1,2),(2,3)\}$. It is 
easy to see that $I(V) = \langle g_1,g_2,g_3\rangle$, where 
$$
g_1(x,y) = x(x-1)(x-2), \quad g_2(x,y) = x(x+1-y), \quad g_3(x,y) = y(x+1-y).
$$
If we use the graded reverse lexicographical order $(n-k-1,k+1)\succ (n-k,k)$,
then $\langle \LT(I) \rangle = \langle y^2, xy, x^3\rangle$ and 
$\RR[V] = \sspan\{1, x, y, x^2\}$. If we use the graded lexicographical order 
$(n-k,k) \succ (n-k-1,k+1)$, then $\langle \LT(I) \rangle = \langle x^2, xy, 
y^3\rangle$ since we also have $I(V) = \langle g_2,g_3, g_4\rangle$ where 
$g_4(x,y) = y^3-6y^2+5y+6x$, and  $\RR[V] = \sspan\{1, x, y, y^2\}$. \qed

\medskip

For $d =2$, the property \eqref{eq:2.1} of $\Lambda$ shows that the set 
$\Lambda$ must be of a stair shape as the lattice points in the unshaded 
area depicted in Fig. 1.

\medskip
\centerline{
\epsfxsize= 2in
\epsffile{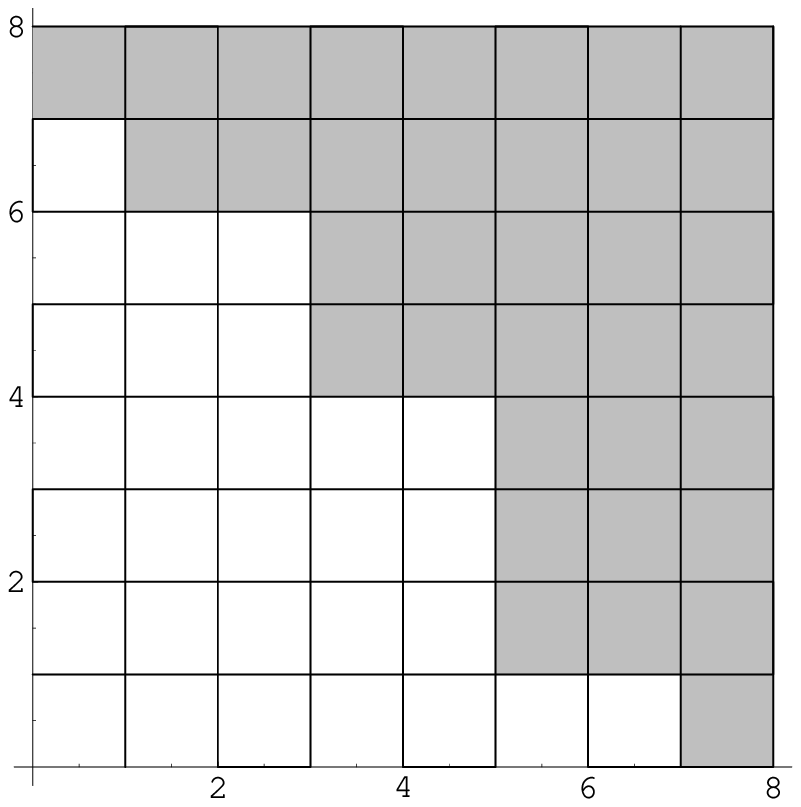}
\hskip .5in
\epsfxsize= 2in
\epsffile{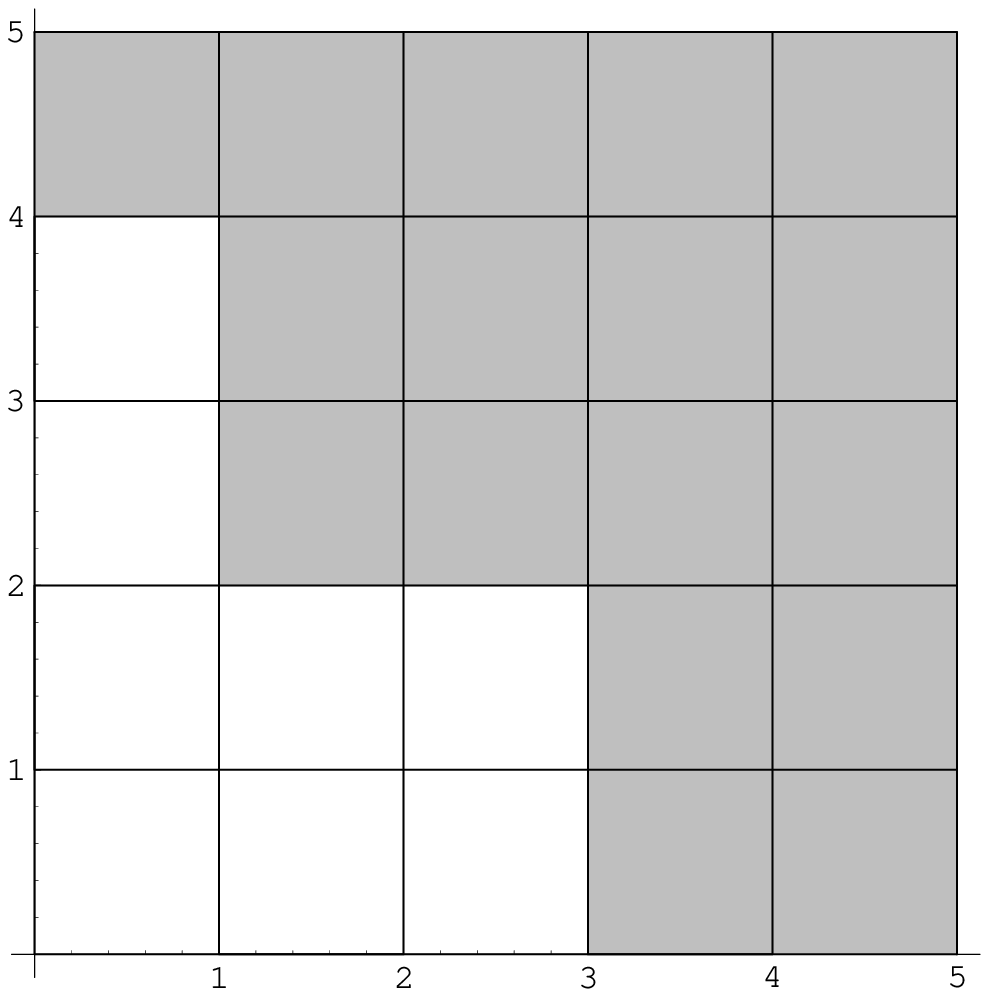} }
\centerline{ Figure 1}
\medskip

\noindent
That is, in the case of $d =2$, for each set $\Lambda$ there is a sequence of
positive integers $n_i$, which satisfies $n_m \le n_{m-1}\le \ldots\le n_0$ 
(some of the $n_i$ can be positive infinity), such that 
\begin{equation} \label{eq:2.3}
 \Lambda = \{(k,l): 0 \le l \le m, 0 \le k \le n_l\}.  
\end{equation} 

\medskip\noindent
{\bf Example 2.2}. Let $\Lambda$ be the lattice set in the two figures in 
Figure 1. For the left one, $m=6$ and 
$\Lambda = \{(i,j),0\le i\le n_j, 0 \le j \le 6\}$ 
 with $(n_0,\ldots,n_6) = (6,4,4,4,2,2,0)$. For the right one, $m = 3$ and
$$
   \Lambda = \{(i,0):0 \le i \le 2\} \cup\{(i,1):0 \le i \le 2\} \cup
     \{(0,2)\} \cup \{(0,3)\} 
$$
with $(n_0,n_1,n_2,n_3) = (2,2,0,0)$. \qed

\medskip

\begin{proposition}
There exists a point set $V$ for which $\RR[V]$ in \eqref{eq:2.2} is given by
the index set $\Lambda$ in \eqref{eq:2.3}.
\end{proposition}

\begin{proof}
Let $x_0, x_1, \ldots, x_{n_0}$ and $y_0, y_1,\ldots,y_m$ be isolated real
numbers. We define the point set $V$ as follows:
$$ 
  \begin{matrix} V = \{(x_0,y_0), & (x_1,y_0),& \ldots,&
      (x_{n_0},y_0),  \\ 
  \,\;\;\,\qquad (x_0,y_1), & (x_1,y_1),& \ldots,& (x_{n_1},y_1) \\ 
   \,\;\;\,\qquad  \vdots & \vdots & \ldots & \vdots&&\\
    \,\;\;\,\qquad   (x_0,y_m), & (x_1,y_m),&\ldots,& (x_{n_m},y_m)\}. 
\end{matrix} 
$$
If $n_0$ is finite, then $|V| = n_0 + n_1 + \ldots + n_m + m+1$. We let 
$n_{m+1} =-1$ and adopt the convention that $\prod_{i=0}^{-1}a_i =1$. Define 
polynomials 
$$
g_k(x) = \prod_{i=0}^{n_k} (x-x_i)\prod_{j=0}^{k-1}(y-y_j), 
\qquad 0\le k \le m+1,
$$
where if $n_k = \infty$, then we take $g_k(x) =1$. Then it is easy to see 
that $\langle I \rangle = \langle g_0, g_1,\ldots, g_{m+1} \rangle$, so that 
$$
\langle \LT(I) \rangle = \langle x^{n_0+1}, x^{n_1+1}y,\ldots,
    x^{n_m+1} y^{m}, y^{m+1}\rangle,  
$$
from which it follows that $\RR[V]= \{x^k y^l: (k,l)\in \Lambda\}$ with
$\Lambda$ given in \eqref{eq:2.3}.
\end{proof}

\medskip\noindent
{\bf Example 2.3}. Consider the ``triangle'' point set 
$$
 \begin{matrix} V = \{(x_0,y_0), & (x_0,y_1),&\ldots, & (x_0,y_m); \\ 
                   & (x_1,y_1) & \ldots,& (x_1,y_m)\\
                   & &  \ddots & \vdots \\
                  & & &  (x_m,y_m)\},  \end{matrix}
$$ 
where all points are isolated and $|V|=(m+1)(m+2)/2$. In this case, 
$n_l = m - l$ for $0 \le l \le m$ 
(making a transpose of the array as in matrix transpose), so that 
$\Lambda = \{(k,l): 0\le k + l \le m\}$ and $\RR[V]=\Pi_m^2$. \qed 

\medskip

\medskip\noindent
{\bf Example 2.4}. Consider the ``product'' point set
$$
V = X \times Y  = \{(x_i,y_j): 0 \le i \le n, 0 \le j \le m\},
$$
for which $|V| = (n+1) (m+1)$. In this case $n_0 = n_1 = \ldots =n_m =n$ and 
$\RR[V] = \{x^k y^l: 0 \le k \le n, 0 \le l \le m\} = \Pi_n^1 \times
   \Pi_m^1$.  \qed

\medskip

The above discussion can be extended to $\RR[x_1,\ldots,x_d]$. We stick to 
the case $d =2$ to keep the notation simple. 

A couple of remarks are in order. First, if $m = \infty$ in Example 2.2, then
$\RR[V] = \RR[x,y]$. This is also the case for Example 2.3 when both $n$ and 
$m$ are infinity. If, however, $n$ is infinite and $m$ is finite (or other way
round), then $\Lambda = \{(k,l): 0 \le l \le m, k\ge 0\}$ is an infinite set
but not all $\NN_0^2$ and $\RR[V] = \Pi^d \times \Pi_m$. Furthermore, the
space of polynomials of one variable appears as a special case.

\medskip\noindent
{\bf Example 2.5}. If $V$ is a point set on the $x$ coordinate line, 
$V = \{(x_0,0), \ldots, (x_n,0)\}$ ($n$ can be infinity), then 
$\RR[V] = \Pi_n^1$. \qed 

\medskip

Second, we should point out that the point set $V$ in the proposition, 
and the above examples, are the simplest examples for which $\RR[V]$ can be
determined. For a generic point set $V$, the problem of determining $\RR[V]$ 
is highly non-trivial. One possible algorithm, at least for $V$ is finite and 
moderate in size, is to check the rank of the matrices whose rows are the 
vectors $X_\alpha:=\{x^\alpha: x \in V\}$. 

The algorithm goes as follows: Fix a graded monomial order, 
starting with $\alpha \in \NN_m:=\{\beta: |\beta|\le m\}$ for a small $m$, so 
that the resulting matrix $(X_\alpha)_{\alpha \in \NN_m}$ has full rank. Then 
add new rows $X_\alpha$ according to the order in $|\alpha| = m+1$ to the
matrix. For each new row added, check the rank of the new matrix; if it has 
full rank, add the next $X_\alpha$ and proceed; if it does not have full rank, 
remove this row and add the next $X_\alpha$ and continue. When the matrix 
becomes a square nonsingular matrix, the corresponding set of $x^\alpha$ will 
be a basis for $\RR[V]$. 

\section{Discrete orthogonal polynomials}
\setcounter{equation}{0}

Let $V$ be a set of isolated points in $\RR^d$. Let $W$ be a real function 
on $V$ and $W(x) \ne 0$ for any $x \in V$. Assume that 
$$
   \sum_{x \in V} |x^\alpha| |W(x)| < \infty \qquad \hbox{for all 
      $\alpha \in \NN_0^d$}
$$
in the case where $V$ is an infinite set. We define a bilinear form 
$\langle \cdot,\cdot \rangle$ on $\Pi^d \times \Pi^d$ by 
$$
\langle f,g \rangle = \CL (f g), \quad \hbox{where} \quad
\CL(f):= \sum_{x \in V} f(x)W(x). 
$$ 
If $\langle f,g \rangle =0$, we say that $f$ and $g$ are orthogonal to each 
other with respect to $W$ on the discrete set $V$. The notation $\CL(f g)$ is 
more convenient for the matrix operations below.

Fix a graded monomial order. Let $\RR[V]$ and $\Lambda = \Lambda(V)$ be defined
as in the previous section (see \eqref{eq:2.2}). Let $\nb = \max \{|\alpha|: 
\alpha \in \Lambda(V)\}$ and $\Lambda_k(V) = \{\alpha \in \Lambda(V):
|\alpha|=k\}$, $0 \le k \le \nb$. Note that $\nb$ can be infinity. Define 
$r_k = |\Lambda_k(V)|$. We denote by $\xb^\Lambda$, $\xb^\Lambda_k$ and 
$\xb^k$, $0 \le k \le \nb$, the sets 
$$
\xb^\Lambda = \{x^\alpha: \alpha \in \Lambda(V)\}, \quad  
\xb_k^\Lambda = \{x^\alpha \in \xb^\Lambda: |\alpha| \le k\}, \quad  
\xb^k = \{x^\alpha: \alpha \in \Lambda_k(V)\}, 
$$  
respectively. We also regard them as column vectors in which the elements are
ordered according to the fixed graded monomial order. 

The set of orthogonal
polynomials on $V$ will be denoted by $\{P_\alpha$, $\alpha \in \Lambda(V)\}$, 
where $P_\alpha$ has degree $|\alpha|$. We introduce the following notion. 
If $\{P_\alpha: \alpha \in \Lambda(V)\}$ is a sequence of polynomials in 
$\RR[V]$, then set $\PP_k:= \{P_\alpha: \alpha \in \Lambda_k(V)\}$. Just as 
in the case of $\xb^k$, we also regard $\PP_k$ as a column vector.  

\begin{definition} \label{def:3.1}
Let $V$ be the set of isolated points and $W$ be a nonzero real function on 
$V$ as above. A sequence of polynomials $\{P_\alpha
\in \Pi_{|\alpha|^d}: \alpha \in \Lambda(V)\}$ is orthogonal with respect 
to $W$ if 
$$
 \CL(\xb^l\PP_k^T)= \langle \xb^l, \PP_k\rangle = 0, \quad k > l, \quad 
\hbox{and}\quad
  \CL(\xb^k \PP_k^T) =\langle \xb^k, \PP_k\rangle = S_k,
$$
for $0 \le k \le \nb$, where $S_k$ is an invertible matrix of size 
$r_k\times r_k$. The sequence is orthonormal with respect to $W$ on $V$ 
if $\langle \PP_k,\PP_k\rangle = I_{r_k}$, the identity matrix, 
for $0 \le k \le \nb$.
\end{definition}

The notation $\CL(\xb^l\PP_k^T)$ is more convenient than $\langle \xb^k,\PP_k
\rangle$, since it shows clearly that this is a matrix of size $r_k\times r_k$.
The orthogonality of $P_\alpha$ is defined as orthogonal to lower degree 
polynomials, as in the continuous case. The polynomials of the same degree 
may not be pairwise orthogonal. By definition, we can write 
$$
  \PP_k = G_k \xb^k + G_{k-1,k} \xb^{k-1} + \ldots + G_{1,k} \xb^0, 
$$
where $G_k$ is a $r_k \times r_k$ matrix, called the leading coefficient
of $\PP_k$. Assume a sequence of orthogonal polynomials $P_\alpha$
exists on $V$. Then we can follow the proof in \cite[Section 3.1]{DX} to 
show that $\{\PP_0,\ldots,\PP_k\}$ is a basis for $\RR[V] \cap \Pi_k^d$, 
$H_k: = \CL(\PP_k \PP_k^T)$ and $G_k$ are both invertible. Furthermore,
the following theorem still holds. 

\begin{theorem} \label{thm:3.2}
A sequence of orthogonal polynomials $\{P_\alpha: \alpha \in \Lambda(V)\}$ 
with respect to $W$ on $V$ exists if and only if the matrices $M_k: = 
\langle \xb_k^\Lambda, \xb_k^\Lambda\rangle$ are nonsingular for 
$0 \le k \le \nb$. 
\end{theorem}

The proof is based on linear algebra and follows exactly as in the continuous 
case; see \cite[Theorem 3.1.6]{DX}. 

If $W$ is positive on $V$, $W(x) >0$ for all $x \in V$, then 
$\langle f, f\rangle =0$ implies $f \equiv 0$ on $V$, so that the bilinear 
form $\langle \cdot, \cdot\rangle =0$ becomes an inner product on $\RR[V]$. 
For such a $W$, orthogonal polynomials on $V$ exist. Furthermore, in this
case, we can have orthonormal bases. 

\begin{theorem} \label{thm:3.3}
If $W$ is a positive function on $V$, then a sequence of orthonormal 
polynomials $\{P_\alpha: \alpha \in \Lambda(V)\}$ with respect to $W$ on 
$V$ exists.
\end{theorem}

\begin{proof}
In this case, the matrix $M_k$ is positive definite since for any nonzero 
column vector $\cb$, $\cb^T M_k \cb = \langle \cb^T \xb^k,\cb^T \xb^k 
\rangle>0$. In particular, $M_\Lambda:=\langle \xb^\Lambda, \xb^\Lambda\rangle$
is a symmetric and positive definite matrix. 
It follows that it can be factored as $M_\Lambda = S D S^T$ where $S$
is a nonsingular lower triangular matrix and $D = \diag\{d_1,d_2,\ldots, 
d_{|V|}\}$ with all $d_i >0$. Let $D^{-1/2} = \diag\{d_1^{-1/2},\ldots,
d_{|V|}^{-1/2}\}$ and $R = D^{-1/2} S^{-1}$. Then $R M_\Lambda R^T = 
\langle R \xb^\Lambda, R \xb^\Lambda\rangle = I_{|V|}$. Since $S$ is 
lower triangular and $\xb^\Lambda = \{\xb^0, \xb^1,\ldots, \xb^\nb\}$ as a 
column vector, we can write the components of $R \xb^\Lambda$ as $\PP_0,
\ldots, \PP_\nb$ where $\PP_k$ consists of polynomials of degree $k$. These
are the orthonormal polynomials. 
\end{proof}

The proof of the theorem provides an algorithm that can be used to construct 
discrete polynomials in several variables. For $|V|$ is finite and moderate in
size, it is rather effective; an example is given in the following section
(Example 4.3).

Let $\CV_k(W)$ denote the space of orthogonal polynomials of degree $k$; that 
is, $\CV_k(W) = \sspan \PP_k$. Evidently $\dim \CV_k(W) = r_k$. Comparing to 
the orthogonal polynomials in the continuous case, the numbers $r_k$ depend
on $V$ and there is no closed formula for them. Furthermore, $r_k \le 
\dim \CP_k^d$ and the equality often does not hold. In fact, $r_k$ may no 
longer be an increasing sequence. 

\medskip\noindent
{\bf Example 3.1}. Let $V =X\times Y$ be the product point set in Example 2.3,
where $X=\{x_0,\ldots,x_n\}$ and $Y = \{y_0, \ldots, y_m\}$. Here $\nb = 
n+m$. Let $W(x,y) = w_1(x) w_2(y)$, where $w_1$ is positive on $X$ and 
$w_2$ is positive on $Y$. Then the orthogonal polynomials $\PP_k$, 
$0 \le k \le \nb$ exist and can be constructed as follows: Let 
$\{p_0, p_1, \ldots, p_n\}$ and $\{q_0, q_1, \ldots, q_n\}$ denote the 
discrete orthogonal polynomials with respect to $w_1$ on $X$ and $w_2$ on 
$Y$, respectively. Then the orthogonal polynomials on $V$ are given by 
$P_{k,l}(x,y)= p_k(x) q_l(y)$. We assume $m$ is finite and $n > m$. Then in 
our vector notation,
$$
  \PP_k(x,y) = \{p_0(x) q_k(y), p_1(x)q_{k-1}(y), \ldots, p_k(x)q_0(y)\} 
\qquad \hbox{for $0 \le k \le m$}, 
$$ 
so that $r_k = k+1$ for $0 \le k \le m$,  
$$
\PP_k(x,y) = \{p_{k-m}(x) q_m(y), \ldots, p_k(x)q_0(y)\} 
\qquad \hbox{for $m+1 \le k \le n$}, 
$$ 
so that $r_k =m+1$ for $m+1 \le k \le n$, 
and
$$
 \PP_k(x,y) = \{p_{k-m}(x) q_m(y), \ldots, p_n(x)
         q_{k-n}(y)\} \qquad \hbox{for $n+1 \le k \le n+m$}, 
$$
so that $r_k =n+m-k+1$ for $n+1 \le k \le n+m-1$. \qed 

\medskip

Next we consider the three-term relations satisfied by the orthogonal 
polynomials. If $P_\alpha$ is an orthogonal polynomial, then $P_\alpha \in
\RR[V]$ so that it is a linear combination of $x^\alpha$ for $\alpha \in 
\Lambda(V)$. Clearly, multiplying by a coordinate $x_i$ gives a polynomial 
$x_i P_\alpha$ of degree $|\alpha|+1$. However, unlike the continuous case,
$x_i x^\alpha$ may not belong to $\RR[V]$ for some $\alpha \in \Lambda(V)$. 
Nevertheless, it is congruent modulus $I(V)$ to a unique polynomial in 
$\RR[V]$. Recall that $\nb = \max \{|\alpha|: \alpha \in \Lambda(V)\}$.

\begin{proposition} \label{prop:3.4}
Let $I : = I(V)$. 
For $0 \le k \le \nb-1$, there exist matrices $A_{k,i}:r_k \times r_{k+1}$,
$B_{k,i}:r_k \times r_k$, and $C_{k,i}:r_k \times r_{k-1}$, such that
for $1 \le i \le d$,
\begin{equation} \label{eq:3.1}
x_i \PP_k(x) = A_{k,i} \PP_{k+1}(x) + B_{k,i} \PP_k(x) + C_{k,i}\PP_{k-1}(x)
\mod I, 
\end{equation}
where $0 \le k \le \nb-1$ and we define $\PP_{-1} = 0$, $A_{k,i} =0$ and 
$C_{-1,i} =0$; moreover, 
\begin{equation} \label{eq:3.2}
A_{k,i}H_{k+1} =\CL( x_i \PP_k \PP_{k+1}^T) =  H_k C_{k+1,i}^T, 
\quad
B_{k,i}H_k = \CL( x_i \PP_k \PP_k^T).
\end{equation} 
\end{proposition}

\begin{proof}
If all components of $x_i \PP_k$ are in $\RR[V]$, this is proved as in the
usual case, by writing $x_i \PP_k$ in terms of orthogonal polynomials  
$\PP_0, \PP_1, \ldots, \PP_{k+1}$ and then use orthogonality. If $P_\alpha
\in \RR[V]$ but $x_i P_\alpha \notin \RR[V]$ for some $|\alpha|=k$, then there
exist a $Q$ in $I(V)$ and an $R_\alpha \in \RR[V]$ such that $x_i P_\alpha(x)=
Q(x) + R_\alpha(x)$, and the degree of $R_\alpha$ is at most $k+1$. We can 
write $R_\alpha = \ab_{k+1} \PP_{k+1} + \ab_k \PP_k + \ab_{k-1} \PP_{k-1}+
\ldots$, where $\ab_j$ are some row vectors of appropriate size. By the 
orthogonality, we get
$$
\ab_{k+1} \CL( \PP_{k+1}\PP_{k+1}^T)
  =\CL( R_\alpha\PP_{k+1}^T ) = \CL( P_\alpha \PP_{k+1}^T),
$$
since $Q$ vanishes on $V$ and $\langle Q, P \rangle = 0$ for any $P$. 
Similarly, we get $\ab_k H_k = \CL( P_\alpha \PP_k^T)$, $\ab_{k-1} H_{k -1}=
\CL( P_\alpha \PP_{k-1}^T)$, and all other $\ab_j$ are equal to zero. In 
vector and matrix notation, this is the three-term relation. The presence of
$Q$ means that the equality holds under modulus $I(V)$ in general.
\end{proof}

\begin{corollary}
Let $I : = I(V)$. If $\{\PP_k\}$ are orthonormal polynomials, then for 
$1 \le i \le d,$
\begin{equation} \label{eq:3.3}
x_i \PP_k(x) = A_{k,i} \PP_{k+1}(x) + B_{k,i} \PP_k(x) 
    + A_{k-1,i}^T \PP_{k-1}(x) \mod I, 
\end{equation}
where $0 \le k \le \nb-1$, $\PP_{-1} = 0$, $A_{k,i} =0$ and $C_{-1,i} =0$; 
moreover, $B_{n,i}$ are symmetric. 
\end{corollary}

In the case of continuous orthogonal polynomials, the matrix $A_{n,i}$ has
more columns than rows and it has full rank. This is no longer true in the 
discrete case. Since $r_k$ may no longer be an increasing sequence, the
matrix $A_{n,i}$ can have more rows than columns; moreover, it may not have
full rank.  

\medskip\noindent
{\bf Example 3.1$^*$}. We continue the example $V =X\times Y$ in Example 3.1.
Assume that $\{p_0,p_1,\ldots,p_n\}$ satisfies the three-term relation 
$$
 x p_k(x) = a_n p_{k+1}(x) + b_k p_k(x) + c_k p_{k-1}(x)
$$
and $\{q_0,q_1,\ldots,q_m\}$ satisfies the three-term relation 
$$ 
 y q_k(y) = a_n' q_{k+1}(y) + b_k' q_k(y) + c_k' q_{k-1}(y), 
$$
respectively. With $x = x_1$ and $y=x_2$, the matrices $A_{k,1}$ and $A_{k,2}$ 
take the form
$$
 A_{k,1} =\left[\begin{matrix} 0&a_0&&\bigcirc\cr
       \vdots&&\ddots& \cr  0&\bigcirc&&a_k \end{matrix} \right]
\quad \hbox{and}\quad
 A_{k,2} =\left[ \begin{matrix} a_k'&&\bigcirc&0\cr &\ddots&&\vdots \cr
          \bigcirc&&a_0'&0\end{matrix} \right],
  \quad \hbox{$0 \le k < m$},
$$
of size $(k+1)\times(k+2)$; 
$$
 A_{k,1} =\left[\begin{matrix} a_{k-m}&&\bigcirc \cr
         &\ddots& \cr \bigcirc&&a_k \end{matrix} \right]
\,\, \hbox{and}\,\,
 A_{k,2} =\left[ \begin{matrix} 0 &0&\ldots & 0\\
         a_{m-1}'&0 & & \bigcirc \\
            &\ddots& \ddots& \cr
          \bigcirc& \ldots&a_0'&0\end{matrix} \right],
  \quad \hbox{$m \le k <n $},
$$
of size $(m+1)\times(m+1)$; and 
$$
 A_{k,1} =\left[\begin{matrix} a_{k-m}&&\bigcirc \cr
         &\ddots& \\ \bigcirc&&a_k \\
         0& \ldots & 0 \end{matrix} \right]
\,\, \hbox{and} \,\,
 A_{k,2} =\left[ \begin{matrix} 0 &\ldots & 0\\
         a_{m-1}'& & \bigcirc \\
            &\ddots&  \cr
          \bigcirc& & a_{k-n}'\end{matrix} \right],
  \quad \hbox{$n \le k <n+m $},
$$
of size $(n+m-k+1)\times (n+m-k)$. Note that $A_{k,2}$ in the case 
$m \le k < n$ does not have full rank. For $n \le k < n+m$, $A_{k,i}$ 
has more rows than columns.
 \qed

\medskip 
One consequence of the three-term relation is the Christoffel-Darboux formula,
$$
  \sum_{j=0}^k \PP_k^T(x) \PP_k(y) = \frac{\PP_{k+1}(x) A_{k,i}^T \PP_k(y)
      - \PP_k(x) A_{k,i} \PP_{k+1}(y)}{x_i - y_i}, \quad \mod I(V),
$$
where $1 \le i \le d$ and $0\le k < \nb$. The proof follows as in the 
continuous case. 

The composite matrix $A_k = (A_{k,1}^T, \ldots, A_{k,d}^T)^T$ plays an 
important
role in Favard's theorem of several variables. This matrix is of size 
$d r_k \times r_{k+1}$. 

\begin{proposition} \label{prop:3.6}
For $0 \le k \le \nb-1$, $d r_k \ge r_{k+1}$; the composite matrix $A_k$ of 
$A_{k,1}, \ldots, A_{k,d}$ and the composite matrix $C_{k+1}$ of $C_{k+1,1}, 
\ldots, C_{k+1,d}$ both have full rank,
\begin{equation} \label{eq:3.4}
  \rank A_k = \rank C_{k+1}^T = r_{k+1}.
\end{equation} 
\end{proposition}

\begin{proof}
Recall that $G_k$ denotes the leading coefficient of $\PP_k$ and it is an 
invertible matrix. Let $e_i = (0,\ldots,0,1,0,\ldots,0)$ be the $i$-th element
of the standard Euclidean basis. By \eqref{eq:2.1}, $\alpha \in \Lambda(V)$ 
implies that $\alpha - e_i \in \Lambda(V)$. If $\alpha \in \Lambda_k(V)$ but 
$\alpha+e_i \notin \Lambda_{k+1}(V)$, then $x_i x^\alpha = x^{\alpha+e_i}\in
\langle \LT(I)\rangle$. We then define the matrix $L_{k,i}$ by
$$
   x_i \xb^k = L_{k,i} \xb^{k+1} \mod \langle \LT(I) \rangle, 
     \qquad 1 \le i \le d. 
$$
The matrix $L_{k,i}$ is of the size $r_k \times r_{k+1}$ and it is uniquely
determined. Comparing the coefficients of $\xb^{k+1}$ in both sides of the 
three-term relation, we see that
\begin{equation} \label{eq:3.5}
  G_k L_{k,i} = A_{k,i} G_{k+1}, \qquad 1 \le i \le d.
\end{equation}
The entries of the matrix $L_{k,i}$ are mostly 0 with at most one 1 in each
row. However, $L_{k,i}$ may not have full rank. 

Since $\alpha \in \Lambda_{k+1}(V)$ implies that $\alpha - e_i \in \Lambda_k$ 
whenever $\alpha_i - 1 \ge 0$, it follows that $d r_k \ge r_{k+1}$ since 
$r_k = |\Lambda_k(V)|$. Moreover, since the column vector $(x_1 \xb^k,\ldots,
x_d \xb^k)$ is equal to $L_k \xb^{k+1}$, where $L_k$ is the composite matrix 
of $L_{k,1},\ldots,L_{k,d}$, and clearly 
$\{x_i \xb^k: 1 \le i \le d\}$ includes $\xb^{k+1}$ as a subset, it follows
that $L_k$ has full rank $r_{k+1}$. 

The equation \eqref{eq:3.5} implies that $A_kG_{k+1} = \diag\{G_k,\ldots, 
G_k\} L_k$. Since $G_k$ invertible implies $\diag\{G_k, \ldots, G_k\}$ 
invertible, it follows that $\rank A_{k+1} = r_{k+1}$. Furthermore, 
\eqref{eq:3.2} implies that $A_k H_{k+1} = \diag\{H_k, \ldots, H_k\} 
C^T_{k+1}$ and $H_k$ is invertible; hence, $\rank C_{k+1} = \rank A_k$.
\end{proof}

Since the matrix $A_k$ has full rank and $d r_k \ge r_{k+1}$, it has a 
generalized 
inverse, $D_k^T$, which is of the size $r_{k+1} \times d r_k$ and can be 
assumed to be of the form $D_k^T = (D_{k,1}^T, \ldots, D_{k,d}^T)$, where 
$D_{k,i}^T$ are of the size $r_{k+1} \times r_k$. Then 
$$
  D_k^T A_k = D_{k,1}^T A_{k,1} + \ldots +  D_{k,d}^T A_{k,d} = I_{r_{k+1}}.  
$$
We note that the generalized inverse is in general not unique. Using $D_n^T$, 
we get from the three-term relation a recursive formula
\begin{equation} \label{eq:3.6}
\PP_{k+1} = \sum_{i=1}^d x_i A_{k,i} \PP_k - \sum_{i=1}^d  B_{k,i} \PP_k -
 \sum_{i=1}^d  C_{k,i} \PP_{k-1},  
\end{equation}
which allows us  to compute $\PP_{k+1}$ using $\PP_k$ and $\PP_{k-1}$. This 
formula is useful in the proof of the analog of Favard's theorem.

According to Propositions 3.4 and 3.6, orthogonal polynomials on a set $V$ 
of isolated points satisfy a three-term relation whose coefficient satisfies
a rank condition. We want to establish that the converse is also true, that
is, an analog of Favard's theorem. To this end, we start with a sequence of
polynomials that satisfies the three-term relation \eqref{eq:3.1} and the 
rank condition \eqref{eq:3.4}, and show that there exist a set $V$ of 
isolated points and a weight function $W$ on $V$ with respect to which 
the polynomials are orthogonal. 

For this purpose let us start with an ideal $I \subset \Pi^d$ and let 
$\Lambda : = \{\alpha \in \NN_0^d: x^\alpha \notin \langle \LT(I) \rangle\}$.  
The proof of Proposition \ref{prop:2.2} shows that $\Lambda$ satisfies the 
property \eqref{eq:2.1}. Assume that there is a sequence of polynomials 
$P_\alpha \in \CS_I: = \sspan\{x^\alpha| x^\alpha \notin \langle \LT(I)
\rangle\}$, where $P_\alpha$ is indexed by $\alpha  
\in \Lambda$ such that $P_0 =1$ and $P_\alpha \in \Pi_{|\alpha|}^d$. 
Set $\nb = \max\{|\alpha|: \alpha \in \Lambda\}$, which can be infinity. 
Let $\Lambda_k = \{\alpha \in \Lambda: |\alpha|=k\}$ and let $\PP_k =
\{P_\alpha: \alpha \in \Lambda_k\}$ for $0 \le k \le \nb$, and 
regard $\PP_k$ as column vectors according to a fixed graded monomial 
order. 

\begin{theorem}
Let $I$ be an ideal of $\Pi^d$ and let $\Lambda$ and $\PP_k$ be as above. 
Assume that $\PP_k$ satisfies the 
three-term relation \eqref{eq:3.1} whose coefficient matrices satisfy the 
rank condition \eqref{eq:3.4}. 
\item{(i)} There is a linear functional $\CL$ on $\CS_I$ for which 
$P_\alpha$ are orthogonal polynomials with respect to the bilinear 
form $\CL(f g) = \langle f, g\rangle$. 
\item{(ii)} If $\nb$ is finite then there exist a set $V$ of isolated points
and a real function $W$ on $V$, $W(x) \ne 0$ for all $x \in V$, such that 
$\{P_\alpha: \alpha \in \Lambda\}$ is a sequence of orthogonal polynomials 
with respect to $W$ on $V$. 
\end{theorem}
 
\begin{proof}
(i) The proof follows along the line of the proof of Favard's theorem for 
continuous orthogonal polynomials of several variables. We shall be brief 
whenever proofs in the two  cases are essentially the same. Using induction, 
it follows from the three-term relation and the rank condition that the 
leading coefficients $G_k$ of $\PP_k$ are invertible. The linear functional, 
$\CL$, defined by 
$$
  \CL 1 =1 \qquad\hbox{and} \qquad \CL(\PP_k) = 0, \quad 1 \le k \le \nb,
$$ 
is well-defined for $\CS_I$, since every polynomial in $\CS_I$ takes
the form $\sum_{\alpha \in \Lambda} c_\alpha x^\alpha$ and $\Lambda$ 
satisfies \eqref{eq:2.1}. Furthermore, using \eqref{eq:3.6}, one can show
by induction that $\CL$ satisfies $\CL (\PP_k \PP_j^T) = 0$ for $k \ne j$, 
and the matrix $H_k = \CL(\PP_k \PP_k^T)$ is invertible. Consequently, 
$\PP_k$ are orthogonal polynomials with respect to the bilinear form 
$\CL(fg) = \langle f, g\rangle$. 

\noindent
(ii) We only need to show that the linear functional $\CL$ can be represented
by a sum over a set of isolated points; that is, $\CL$ can be written as 
\begin{equation} \label{eq:3.7}
  \CL f = \lambda_1 f(\xb_1) + \ldots + \lambda_N f(\xb_N)
\end{equation}  
for some $\xb_i$ with $\lambda_i \ne 0$ for $1 \le i \le N$. Assume that $\CL$
has such an expression. If $\xb_i$ are known then $\CL P_\alpha = 
\delta_{\alpha,0}$, $\alpha \in \Lambda$, becomes a system of equations on 
$\lambda_i$. The coefficient matrix of this system is $[P_\alpha(\xb_i)]$,
where $\alpha \in \Lambda$ and $1 \le i \le N$. In particular, if $N =  
|\Lambda|$ then the matrix is a square matrix. Its determinant is a polynomial 
in variables $\xb_1,\xb_2,\ldots, \xb_N$ and defines a hypersurface in 
$\RR^{dN}$. Hence, for almost all choices of the values of $\xb_1,\ldots,
\xb_N$, the determinant is non-zero. Furthermore, by Cremer's rule, it is 
possible to choose a set $V = \{\xb_1,\ldots, \xb_N\}$ such that $\lambda_i 
\ne 0$ for $1 \le i \le N$. Hence, $P_\alpha$ are orthogonal with respect 
to $W$ on $V$, where $W$ is defined by $W(\xb_i) = \lambda_i$, 
$1 \le i \le N$.   
\end{proof} 

A couple of remarks are in order. First of all, if $\nb$ is infinity, we 
do not know if $\CL$ can be written as a sum of point evaluations, 
$\CL f = \sum_{x \in V} f(x)$, on a countable set $V$. In fact, 
if $I$ is the trivial ideal $\langle 1\rangle$, then $\Lambda = \NN_0^d$; the 
three-term relation takes the same form as that for continuous orthogonal 
polynomials and the rank condition remains the same. Hence, the three-term 
relation and the rank condition are not enough to give further information
on the linear functional. 

The same phenomenon will happen to the 
case that $\CS_I$ is a product polynomial space, say, $\Pi_m^d \times \Pi^d$, 
for which the orthogonal polynomials are $p_j(x) q_k(y)$, $ 0 \le j \le m$
and $k \ge 0$. No matter if $q_k$ are orthogonal with respect to a linear 
functional defined by an integral or to a linear functional defined by an 
infinite sum, the three-term relation will take the same form and the rank 
condition will also remain the same. 

Secondly, the proof of the theorem only shows that $W$ is nonzero at 
every point of $V$. This is enough if we only deal with orthogonality 
but not orthonormality. The theorem simply states that the rank condition 
\eqref{eq:3.4} and the three-term relation \eqref{eq:3.1} are enough
to ensure orthogonality. The following corollary is about the case of 
orthonormality, where we do get positive weight. The difference is in
the three-term relations \eqref{eq:3.3} vs \eqref{eq:3.1}.

\begin{corollary}
Let $I$ be an ideal of $\Pi^d$ and let $\Lambda$ and $\PP_k$ be as in 
the theorem. Assume that $\nb$ is finite. If $\PP_k$ satisfies the 
three-term relation \eqref{eq:3.3} whose coefficient matrices satisfy 
the rank condition \eqref{eq:3.4}, then there is a set $V$ of isolated 
points and a positive function $W$ on $V$ such that $\{P_\alpha: \alpha 
\in \Lambda\}$ is a sequence of orthonormal polynomials with respect to 
$W$ on $V$. 
\end{corollary}

\begin{proof}
According to the theorem, there is a set $V$ so that $P_\alpha$ are orthogonal
with respect to the bilinear form defined by the linear functional $\CL$ of 
the form \eqref{eq:3.7}. We need to show that $H_k = \CL(\PP_k \PP_k^T)$ is
an identity matrix for $0 \le k \le n$. This can be established by induction.
Since $\PP_0 =1$ and $\CL 1 =1$, we have $H_0 =1$. 
By \eqref{eq:3.2} with $C_{k+1}^T = A_k$, $A_k H_{k+1}=\diag \{H_k,\ldots,
H_k\}A_k$. Assume $H_k = I_{r_k}$. Then $\diag \{H_k,\ldots,H_k\}$ is an 
identity matrix, so that $H_{k+1}$ is an identity matrix by the rank condition.
From the fact that $\CL(\PP_k \PP_k^T) = I$, it follows easily that $\CL$
is a positive definite linear functional, which shows in particular that
$W(x_i)= \lambda_i > 0$ for all $x_i \in V$. 
\end{proof}

The coefficient matrices of the three-term relation \eqref{eq:3.3} can be 
used to define the analog of Jacobi matrices, 
$$
 J_i =  \left[ \begin{matrix} B_{0,i}&A_{0,i}&&&\bigcirc\cr
     A_{0,i}^T&B_{1,i}&A_{1,i}&&\cr  &\ddots&\ddots&\ddots&\cr
     &&A_{\nb-3,i}^T& B_{\nb-2,i}&A_{\nb-2,i}\cr
     \bigcirc&&&A_{\nb-2,i}^T&B_{\nb-1,i}\end{matrix} \right],
\qquad 1\le i\le d.
$$
According to the three-term relation, $J_i$ is the matrix representation of 
the operator $x_i : P\mapsto x_i P$. Since $x_ix_j \PP_k = x_jx_i \PP_k 
\mod I(V)$, $1 \le k \le \nb$, it follows that these matrices commute. 
In other words, we have the following proposition.

\begin{proposition}
If $\nb$ is finite, then the Jacobi matrices commute, $J_i J_j = J_j J_i$, 
$1 \le i,j\le d$. 
\end{proposition} 

If $\nb$ is infinity, then the matrices become infinite and we need to consider
them as operators. Still, the Jacobi matrices formally commute. In the case
that $\Lambda = \NN_0^d$, see \cite{X94}. 

\section{Examples}
\setcounter{equation}{0}

Some examples have been given in the previous sections. To make them more 
concrete, let us mention two classical discrete polynomials. See, for
example, \cite{KM1,KS}.   

The Hahn polynomials, $Q(x; a,b,N)$, are discrete orthogonal polynomials 
defined on the set $V=\{0,1,\ldots, N\}$ and are orthogonal with respect to 
the hypergeometric distribution $(a+1)_x (b+1)_{N-x}/(x! (N-x)!)$,  
\begin{align*} 
 \sum_{x=0}^N \binom{x+a}{x}\binom{N-y+b}{N-y} 
   &   Q_n(x; a, b, N)Q_m(x; a, b, N) \\
 = &\, \frac{(-1)^n n! (b+1)_n (n+a+b+1)_{N+1}}
     {N! (2n+a+b+1) (-N)_n (a+1)_n} \delta_{n,m}, \qquad n,m\le N. 
\end{align*}
Their explicit formulas are given in terms of ${}_3F_2$ series,
\begin{equation*}
  Q_n(x; a, b, N) := {}_3 \wt F_2 \left( \begin{matrix} -n, n+a + b+1, -x\\
       a+1, -N \end{matrix}; 1\right), \qquad n = 0, 1, \ldots, N,    
\end{equation*}
where ${}_3 \wt F_2$ is defined as the usual ${}_3F_2$ with the summation 
terminating at $N$. The Hahn polynomials, $Q_n(x) = Q_n(x; a, b, N)$,
satisfy the three-term relation 
$$
  - x Q_n(x) = A_n Q_{n+1}(x) - (A_n+C_n) Q_n(x) +C_n Q_{n-1}(x), 
$$
where 
$$
 A_n = \frac{(n+a+b+1)(n+a+1)(N-n)}{(2n+a+b)(2n+a+b+1)}, \qquad
 C_n = \frac{n(n+b)(n+a+b+N+1)}{(2n+a+b)(2n+a+b+1)}.
$$

The Meixner polynomials, $M_n(x;b,c)$, are discrete orthogonal polynomials 
defined on the set $V=\NN_0$ and are orthogonal with respect to the negative
binomial distribution $(b)_x c^x/x!$, 
$$
  \sum_{x=0}^\infty \frac{(b)_x} {x!} c^x M_m(x;b,c)M_n(x;b,c) =  
       \frac{c^{-n} n!}{(b)_n (1-c)^b} \delta_{m,n}, 
$$
where $(a)_m$ denote the Pochhammer symbol $(a)_m = a(a+1) \ldots (a+m-1)$. 
Their explicit formula is given in terms of ${}_2F_1$ series,
$$
M_n(x):= M_n(x;b,c) = {}_2 \wt F_1(-n ,-x; b; 1- c^{-1}), \qquad 
n =0,1,2,\ldots   
$$
which satisfies the three-term relation
$$
  (c-1)xM_n(x) = c(n+b)M_{n+1} - (n+(n+b)c)M_n(x)+n M_{n-1}(x).
$$

\medskip\noindent
{\bf Example 4.1} As a special case of Example 3.1, we have the product Hahn 
polynomials of two variables, $Q_n(x;a_1,b_1,N)Q_m(x;a_2,b_2,M)$, 
$0 \le n \le N$, $0 \le m\le M$, which are orthogonal on the set 
$V = \{0,1,\ldots,N\} \times \{0,1,\ldots,M\}$ with respect to the weight 
function 
$$
W(x,y) =\binom{x+a_1}{x}\binom{N-x+a_1}{N-x}\binom{y+b}{y}\binom{N-y+b}{N-y},  
$$
and the product of Hahn and Meixner polynomials, $Q_n(x;a_1,a_2,N)M_m(x;b,c)$, 
$0\le n \le N, m \ge 0$, which are orthogonal on the set $V = \{0,1,\ldots,N\} 
\times \NN_0$ with respect to the weight function 
$$
W(x,y) = \binom{x+a_1}{x}\binom{N-x+a_1}{N-x} \cdot
    \frac{(b)_y} {y!} c^y.    
$$
We also have product Meixner polynomials $M_n(x;b_1,c_1)M_m(x;b_2,c_2)$, 
$n \ge 0$, $m \ge 0$, which are orthogonal on the set $V = \NN_0^2$ with 
respect to the weight function $(b_1)_x (b_2)_y c^{x+y}/ (x! y!)$. In the
last case, we have $r_k = k+1$ for all $k \ge 0$, as for the usual continuous
orthogonal polynomials. 

For the product orthogonal polynomials, the coefficient matrices $A_{k,i}$ 
in the three-term relation are given in Example 3.1$^*$. Note that it is 
easy to get orthonormal Hahn and Meixner polynomials (multiply by square 
root of the normalization constant), and the product of the orthonormal 
polynomials gives orthonormal polynomials in two variables. \qed

\medskip\noindent
{\bf Example 4.2}. Let $V$ be the ``triangle'' point set in Example 2.2. Then
the orthogonal polynomials $\PP_k$, $0 \le k \le m$ exist with $r_k = k+1$.  
This is the case that works exactly as in the case of the continuous orthogonal
polynomials. As an example, let us mention the Hahn polynomials of two 
variables as defined in \cite{KM}. They are given by
\begin{align*}
  \phi_{n,m} (x,y;\sigma,N) & = (-1)^{n+m} \frac{(\sigma_1+1)_n (\sigma_2+1)_m 
             (-N+x)_m} {(\sigma_3+1)_m (\sigma_2+\sigma_3+2n+1)_n (-N)_m} \\
      & \times Q_n(x;\sigma_1, \sigma_2+\sigma_3+2n+1, N-m) 
      Q_m(y;\sigma_2,\sigma_3, N-x),  
\end{align*}
and are orthogonal on the set $V = \{(x,y) \in \NN_0^2: 0 \le x+y \le N\}$
with respect to the weight function
$$
 W(x,y)= \binom{x+\sigma_1}{x}\binom{y+\sigma_1}{y} 
   \binom{N-x-y+\sigma_3}{N-x-y}. 
$$
In this case $r_k = k+1$ and the matrix $A_{k,1}$ and $A_{k,2}$ are of the
size $(k+1)\times (k+2)$, just as in the usual continuous case.
\qed

\medskip
For extensions of classical discrete orthogonal polynomials to several 
variables, we refer to \cite{KM,Tr1,Tr2}. One can also extend discrete 
$q$-orthogonal polynomials to several variables. 

In the following we consider an example in which $V$ contains 8 points and
$\Lambda$ is given as in the right figure of Figure 1. 

\medskip\noindent
{\bf Example 4.3}. 
We set   
$$
V =  \{(-1,-1), (0,-1),(1,-1),(-1,0),(0,0),(1,0),(-1,1),(-1,2)\}.
$$ 
Then $\RR[V] = \sspan\{1,x,y,x^2,xy, y^2,x^2y,y^3\}$ as shown in Example 2.2. 
We use the method described in Theorem 3.3 to construct orthogonal polynomials
on $V$ with respect to the linear functional 
$$
\CL(f) = \frac{1}{8}\sum_{x \in V}f(x)
$$  
That is, we compute the matrix $M_\Lambda = \langle x^\Lambda, 
x^\Lambda\rangle$, where $\langle f, g\rangle = \CL(f g)$, and factor it 
as $SDS^T$. The orthogonal polynomials are given as follows:
\begin{align*}
P_0^0(x,y) & = 1,\\
P_0^1(x,y) & = 1+4x,\\
P_1^1(x,y) & = 3+12 x + 22y,\\
P_0^2(x,y) & =  -26+x+35x^2 - 4y, \\
P_1^2(x,y) & = 3 + 3 x + x^2 + 6 y + 8 x y, \\
P_2^2(x,y) & = -20 + 31 x - x^2 + 11 y + 60 x y + 51 y^2, \\ 
P_0^3(x,y) & = -20 + 3 x + 27 x^2 - 45 y + 4 x y + 56 x^2 y - 5 y^2, \\
P_1^3(x,y) & = -9 x + 9 x^2 - 50 y - 12 x y + 12 x^2 y - 30 y^2 + 20 y^3,  
\end{align*} 
where $P_i^k$ is a polynomial of degree $k$. By construction, these polynomials
are mutually orthogonal and become orthonormal upon multiplying by proper 
constants. The corresponding dimensions of $\PP_0,\PP_1,\PP_2,\PP_3$ are 
$1,2,3,2$. 
\qed

\medskip
In the case that $V$ is a set of lattice point and $V$ satisfies 
\eqref{eq:2.1}, we can take $\Lambda  = V$. Example 4.1 and 4.2 are
examples of such a case. If we take $V = \Lambda$ in Example 4.3, we get
orthogonal polynomials $Q_j^k(x,y)=P_j^k(x+1,y+1)$, where $P_j^k$ are those 
in the example. 

\bigskip\noindent
{\it Acknowledgments.} The author would like to thank one referee for 
his careful review and helpful suggestions.

\enddocument